\documentclass[11pt]{article}
\usepackage{latexsym}

\newtheorem{theorem}{Theorem}[section]
\newtheorem{conjecture}[theorem]{Conjecture}

\newtheorem{proposition}[theorem]{Proposition}
\newtheorem{definition}{Definition}
\newtheorem{corollary}[theorem]{Corollary}

\newcommand{\qed}{\ \hfill\mbox{$\Box$}\vspace{\baselineskip}}
\newenvironment{proof}{\noindent {\bf Proof:}}{{\qed}}

\newcommand{\conv}{\mbox{{\rm conv}}}
\newcommand{\vertx}{\mbox{{\rm vert}}}
\begin{document}

\title{A combinatorial study of multiplexes \\ and ordinary polytopes}

\author{Margaret M. Bayer
        \and Aaron M. Bruening\thanks{Current address:
	723 NE Tudor Rd.\ \#4, Lee's Summit, MO  64086}
	\and Joshua Stewart\thanks{Current address:
	5381 Brook Way \#6, Columbia, MD  21044} \\ \\
        Department of Mathematics\\
	University of Kansas\\
	Lawrence, KS  66045-2142}

\date{December 2000}

\maketitle

\begin{abstract}
Bisztriczky defines a multiplex as a generalization of a simplex, and an
ordinary polytope as a generalization of a cyclic polytope.
This paper presents results concerning the combinatorics of multiplexes and
ordinary polytopes.
The flag vector of the multiplex is computed, and shown to equal the flag
vector of a many-folded pyramid over a polygon. 
Multiplexes, but not other ordinary polytopes, are shown to be elementary.
It is shown that all complete subgraphs of the graph of a multiplex determine
faces of the multiplex.
The toric $h$-vectors of the ordinary 5-dimensional polytopes are given.
Graphs of ordinary polytopes are studied.
Their chromatic numbers and diameters are computed, and they are shown to be
Hamiltonian.
\end{abstract}

\section{Introduction}
A {\em convex polytope} is the convex hull of a finite set of points in
${\bf R}^d$.
A $d$-dimensional convex polytope has faces of every dimension from 0
(the vertices) to $d-1$ (the facets).
The set of faces, ordered by inclusion, forms a lattice.
The properties of the face lattice are known as the {\em combinatorial} 
properties of the polytope.
Of particular interest are the number of faces of each dimension (the 
$f$-vector), the graph consisting of the vertices and edges, and the
polyhedral and simplicial complexes that subdivide the polytope.
(The latter are not strictly speaking combinatorial.  Different polytopes
with the same face lattice can have different geometric subdivisions.)

The greatest progress in understanding combinatorial properties of convex
polytopes has been made for the special class of simplicial polytopes.
A {\em simplex} is the convex hull of affinely independent points.
A polytope is {\em simplicial} if all its (proper) faces are simplices.
A milestone was the characterization of the $f$-vectors of simplicial
polytopes, conjectured by McMullen (\cite{McM}), and proved by Billera
and Lee (\cite{billera-lee}) and Stanley (\cite{stanley}).
In this work as well as research on triangulations of simplicial
polytopes, the {\em cyclic} polytopes play a special role.
These are simplicial polytopes, with vertices chosen as points on the
moment curve, $\{(t,t^2, t^3, \ldots, t^d): t\in {\bf R}\}$.
The combinatorial structure of the cyclic polytope depends only on the
number, not the placement, of the points on the moment curve.

The combinatorial study of general  polytopes is hampered by the difficulty
in generating random combinatorial structures.
A polytope in general position, that is, one formed by choosing points at random
and taking their convex hull, is generally simplicial.
An alternative construction, of intersecting a randomly-chosen set of 
hyperplanes,
produces {\em simple} polytopes, which have face lattices dual to simplicial
polytopes.
Various geometric operations (for example, pyramiding, and truncation) can be 
performed on polytopes to produce nonsimplex faces, but the effect of these
on the face lattice is quite controlled.

Bisztriczky (\cite{Bisz}) defines a class of nonsimplicial
polytopes called {\em ordinary} polytopes.
The faces of ordinary polytopes are {\em multiplexes}, which generalize
the simplex in a combinatorial way, but are generally not simplicial.
A special case of the ordinary polytope is the cyclic polytope.
Thus, ordinary polytopes hold out promise of playing an important role in
the combinatorial study of nonsimplicial polytopes.
In this paper we continue the study of multiplexes and ordinary polytopes
begun by Bisztriczky (\cite{Bisz-mult,Bisz}) and Dinh (\cite{Dinh}).

\section{Multiplexes}

Bisztriczky defines a multiplex as a generalization of a 
simplex.
\begin{definition}[\cite{Bisz-mult}]  {\em
A {\em multiplex} is a polytope with an ordered list of vertices,
$x_0$, $x_1$, \ldots, $x_n$, with facets $F_0$, $F_1$, \ldots, $F_n$ given by 
$$F_i = \conv\{x_{i-d+1}, x_{i-d+2}, \ldots, x_{i-1}, x_{i+1}, x_{i+2}, \ldots,
x_{i+d-1}\},$$
with the conventions that $x_i=x_0$ if $i<0$, and $x_i=x_n$ if $i>n$.
}\end{definition}
Note that if $n=d$, then the multiplex is a simplex.
It is easy to check that for $n=d+1$, the multiplex is a $(d-2)$-fold
pyramid over a quadrilateral with vertex set $\{x_0, x_1, x_d, x_{d+1}\}$.
Multiplexes are not all pyramids, however.
Every polygon (two-dimensional polytope) is a multiplex with an appropriate
ordering of vertices.
Following are important results about multiplexes by Bisztriczky.
For $n\ge d\ge 2$,
let $M^{d,n}$ be the $d$-dimensional multiplex with $n+1$ vertices.
\begin{theorem}
[\cite{Bisz-mult}]
\label{multi-backgr}
\begin{enumerate}
\item $M^{d,n}$ exists
      for every $d$ and $n$ with $n\ge d\ge 2$.
\item Every multiplex is self-dual.
\item Every face and every quotient of a multiplex is a multiplex.
\item The number of $i$-dimensional faces of $M^{d,n}$
      is ${d+1\choose i+1}+(n-d){d-1\choose i}$.
\end{enumerate}
\end{theorem}

The $f$-vector of $M^{d,n}$
is the same as the $f$-vector of a certain pyramid over a polygon.
This fact extends to the ``flag vector'' of the multiplex.
A chain of faces, $\emptyset\subset F_1\subset F_2\subset \cdots
\subset F_r\subset P$ is an {\em $S$-flag} if 
$S=\{\dim F_1, \dim F_2, \ldots, \dim f_r\}$.
The number of $S$-flags of a polytope $P$ is written $f_S(P)$, and
the length $2^d$ vector $(f_S(P))_{S\subseteq\{0, 1, \ldots, d-1\}}$
is the {\em flag vector} of $P$.

\begin{theorem}
The multiplex $M^{d,n}$ has the same flag vector as the $(d-2)$-fold
pyramid over the $(n-d+3)$-gon.  
The common flag vector is given by
\begin{eqnarray}
\label{multiplex-flag}
f_S&=&{d+1 \choose \mbox{$s_1+1$, $s_2-s_1$, \ldots, $s_r-s_{r-1}$,
  $d-s_r$}}\\
  & & \times\left[ 1+\frac{n-d}{(d+1)d(d-1)}\sum_{j=1}^r 
  (s_j+1)(s_{j+1}-s_j)(s_{j+1}-1) \right], \nonumber
\end{eqnarray}
where $S=\{s_1,s_2,\ldots, s_r\}$, and $s_{r+1}=d$.
\end{theorem}
\begin{proof}
We first prove that the flag vector of the pyramid is given by 
equation~(\ref{multiplex-flag}).
For $m\ge 0$ let $Q^{d,m}$ be the $(d-2)$-fold pyramid over the $(m+3)$-gon.  
Let $T^d$ be the $d$-simplex.
We prove by induction on $|S|$ that
\begin{eqnarray*}
f_S(Q^{d,m})-f_S(T^d)&=&{d+1 \choose \mbox{$s_1+1$, $s_2-s_1$, \ldots, 
  $s_r-s_{r-1}$, $d-s_r$}}\\
  & & \times\frac{m}{(d+1)d(d-1)}\sum_{j=1}^r 
  (s_j+1)(s_{j+1}-s_j)(s_{j+1}-1).
\end{eqnarray*}
The formula is valid for $S=\emptyset$, since the sum is empty.
We need the $|S|=1$ case as well.
Write $S=\{s_1\}$.
\begin{eqnarray*}
\lefteqn{f_{s_1}(Q^{d,m})-f_{s_1}(T^d) =  
      \sum_{i=-1}^2 {d-2 \choose s_1 - i}f_i(\mbox{$(m+3)$-gon})
		    -{d+1\choose s_1+1}}\\
      & &={d-2\choose s_1+1}+{d-2\choose s_1}(m+3) +{d-2\choose s_1-1}(m+3)
         +{d-2\choose s_1-2}-{d+1\choose s_1+1}\\
      & &= {d-1\choose s_1}m
      ={d+1\choose s_1+1}\frac{m}{(d+1)d(d-1)}(s_1+1)(d-s_1)(d-1).
\end{eqnarray*}
For the induction step we write $f_S$ in terms of $f_{S\setminus\{s_r\}}$.
We need two basic combinatorial facts.  
For a simplex $T^q$, $$f_S(T^q)= 
{q+1 \choose \mbox{$s_1+1$, $s_2-s_1$, \ldots, $s_r-s_{r-1}$,$q-s_r$}}.$$
The $s_r$-faces of the pyramid $Q^{d,m}$ are all either 
$(s_r-2)$-fold pyramids over the $(m+3)$-gon or simplices.
Let ${\cal T}_r$ be the set of $s_r$-faces that are simplices.
Let ${\cal P}_r$ be the set of $s_r$-faces that are pyramids (not simplices); 
the number of such faces is ${d-2\choose s_r-2}$.
Thus,
\begin{eqnarray*}
\lefteqn{
f_S(Q^{d,m})-f_S(T^d)=}\\& & \sum_{F\in {\cal P}_r} f_{S\setminus\{s_r\}}(F)
   +  \sum_{F\in {\cal T}_r} f_{S\setminus\{s_r\}}(F)
   - {d+1 \choose \mbox{$s_1+1$, $s_2-s_1$, \ldots, $s_r-s_{r-1}$,$d-s_r$}}\\
   &=&{d-2\choose s_r-2}f_{S\setminus\{s_r\}}(Q^{s_r,m})\\ & &
      + \left(f_{s_r}(Q^{d,m})-{d-2\choose s_r-2}\right)
      {s_r+1\choose \mbox{$s_1+1$, $s_2-s_1$, \ldots, $s_r-s_{r-1}$}}\\
   & &{} - {d+1 \choose \mbox{$s_1+1$, $s_2-s_1$, \ldots, $s_r-s_{r-1}$,$d-s_r$}}\\
   &=&{d-2\choose s_r-2}(f_{S\setminus\{s_r\}}(Q^{s_r,m})-
   f_{S\setminus\{s_r\}}(T^{s_r}))\\
      & &{} +{s_r+1\choose \mbox{$s_1+1$, $s_2-s_1$, \ldots, $s_r-s_{r-1}$}}
      (f_{s_r}(Q^{d,m})-f_{s_r}(T^d)) 
      \end{eqnarray*} 
      \pagebreak
   \begin{eqnarray*}
   &=&{d-2\choose s_r-2}
      {s_r+1\choose \mbox{$s_1+1$, $s_2-s_1$, \ldots, $s_r-s_{r-1}$}}\\
      & & \times
      \frac{m}{(s_r+1)s_r(s_r-1)}\sum_{j=1}^r(s_j+1)(s_{j+1}-s_j)(s_{j+1}-1)\\
   & & {}  + {s_r+1\choose \mbox{$s_1+1$, $s_2-s_1$, \ldots, $s_r-s_{r-1}$}}
      {d+1\choose s_r+1}\frac{m}{(d+1)d(d-1)}(s_r+1)(d-s_r)(d-1)\\
   &=&{d+1 \choose \mbox{$s_1+1$, $s_2-s_1$, \ldots, 
  $s_r-s_{r-1}$, $d-s_r$}}\\
  & & \times\frac{m}{(d+1)d(d-1)}\sum_{j=1}^r 
  (s_j+1)(s_{j+1}-s_j)(s_{j+1}-1).
\end{eqnarray*}
So the formula for the flag vector of the $(d-2)$-fold pyramid over the
\mbox{$(n-d+3)$-gon} is valid.

Now we show that the same formula holds for the multiplex.  
This is proved by a recursion different from the one used for the pyramid.
Recall the formula for the $f$-vector of the multiplex
$M^{d,n}$ in Theorem~\ref{multi-backgr}.
$$f_i(M^{n,d})= {d+1\choose i+1}+(n-d){d-1\choose i}.$$
This is easily seen to agree with equation~(\ref{multiplex-flag}) for 
$S=\{i\}$.
Next is the proof for the special case where $S=\{0,i\}$.
For $x$ a vertex of $M^{d,n}$, write $[x,M^{d,n}]$ for the quotient of the
polytope $M^{d,n}$ by the vertex $x$.  
This is itself a multiplex of dimension $d-1$.
\begin{eqnarray*}
\lefteqn{
f_{0,i}(M^{d,n}) = \sum_{\mbox{$x$ vertex of $M^{d,n}$}}f_{i-1}([x,M^{d,n}])}\\
   &=&\sum_{\mbox{$x$ vertex of $M^{d,n}$}}{d\choose i}+
   (f_0([x,M^{d,n}])-d){d-2\choose i-1}\\
   &=&(n+1){d\choose i}-(n+1)d{d-2\choose i-1}+
      {d-2\choose i-1}\sum_{\mbox{$x$ vertex of $M^{d,n}$}}f_0([x,M^{d,n}])\\
   &=&(n+1){d\choose i}-(n+1)d{d-2\choose i-1}+
      {d-2\choose i-1}f_{0,1}(M^{d,n})\\
   &=&(n+1){d\choose i}-(n+1)d{d-2\choose i-1}+
      {d-2\choose i-1}2f_1(M^{d,n})\\
   &=&(n+1){d\choose i}-(n+1)d{d-2\choose i-1}+
      {d-2\choose i-1}(d(d+1)+2(n-d)(d-1)\\
   &=&(n+1){d\choose i}+(d-2)(n-d) {d-2\choose i-1}.
\end{eqnarray*}
This can be shown to agree with equation~(\ref{multiplex-flag}).
We now show that equation~(\ref{multiplex-flag}) holds by induction on $|S|$.
\begin{eqnarray*}
\lefteqn{
f_S(M^{d,n})= \sum_{\mbox{$F$ $s_r$-face of $M^{d,n}$}}
           f_{S\setminus\{s_r\}}(F)}\\
	   &=& \sum_{\mbox{$F$ $s_r$-face of $M^{d,n}$}}
           f_{S\setminus\{s_r\}}(M^{s_r,f_0(F)-1}))\\
	   &=& \sum_{\mbox{$F$ $s_r$-face of $M^{d,n}$}}
      {s_r+1\choose \mbox{$s_1+1$, $s_2-s_1$, \ldots, $s_r-s_{r-1}$}}\\
      & & \times\left[1+ \frac{f_0(F)-s_r-1}{(s_r+1)s_r(s_r-1)}
	 \sum_{j=1}^{r-1}(s_j+1)(s_{j+1}-s_j)(s_{j+1}-1)\right]\\
	   &=& {s_r+1\choose \mbox{$s_1+1$, $s_2-s_1$, \ldots, $s_r-s_{r-1}$}}\\
      & & \times\left[1- \frac{s_r+1}{(s_r+1)s_r(s_r-1)}
	 \sum_{j=1}^{r-1}(s_j+1)(s_{j+1}-s_j)(s_{j+1}-1)\right]
	   \sum_{\mbox{$F$ $s_r$-face of $M^{d,n}$}} 1\\
	   &+& {s_r+1\choose \mbox{$s_1+1$, $s_2-s_1$, \ldots, $s_r-s_{r-1}$}}\\
      & & \times\left[\frac{1}{(s_r+1)s_r(s_r-1)}
	 \sum_{j=1}^{r-1}(s_j+1)(s_{j+1}-s_j)(s_{j+1}-1)\right]
	   \sum_{\mbox{$F$ $s_r$-face of $M^{d,n}$}} f_0(F)\\
	   &=& {s_r+1\choose \mbox{$s_1+1$, $s_2-s_1$, \ldots, $s_r-s_{r-1}$}}\\
      & & \times\left[1- \frac{s_r+1}{(s_r+1)s_r(s_r-1)}
	 \sum_{j=1}^{r-1}(s_j+1)(s_{j+1}-s_j)(s_{j+1}-1)\right]
	   f_{s_r}(M^{d,n})\\
	   &+& {s_r+1\choose \mbox{$s_1+1$, $s_2-s_1$, \ldots, $s_r-s_{r-1}$}}\\
      & & \times\left[\frac{1}{(s_r+1)s_r(s_r-1)}
	 \sum_{j=1}^{r-1}(s_j+1)(s_{j+1}-s_j)(s_{j+1}-1)\right] f_{0,s_r}(M^{d,n})\\
	   &=& {s_r+1\choose \mbox{$s_1+1$, $s_2-s_1$, \ldots, $s_r-s_{r-1}$}}\\
      & & \times\left[1- \frac{s_r+1}{(s_r+1)s_r(s_r-1)}
	 \sum_{j=1}^{r-1}(s_j+1)(s_{j+1}-s_j)(s_{j+1}-1)\right]\\
	 & & \times
	   \left[{d+1\choose s_r+1}+(n-d){d-1\choose s_r}\right]\\
	   &+& {s_r+1\choose \mbox{$s_1+1$, $s_2-s_1$, \ldots, $s_r-s_{r-1}$}}
      \frac{1}{(s_r+1)s_r(s_r-1)}\\ & & \times 
	 \sum_{j=1}^{r-1}(s_j+1)(s_{j+1}-s_j)(s_{j+1}-1) 
          \left[(n+1){d\choose s_r}+(d-2)(n-d) {d-2\choose s_r-1}\right].
\end{eqnarray*}
Algebraic manipulation gives equation~(\ref{multiplex-flag}).
\end{proof}

The {\em toric $h$-vector} of a polytope is a length $d+1$ vector whose 
components are linear functions of the flag vector.
When the polytope has rational vertices, the toric $h$-vector is the
sequence of middle perversity intersection homology Betti numbers of
the associated toric variety.
(The multiplex has a rational realization.)
The toric $h$-vector for Eulerian posets (including face lattices of 
polytopes) is defined by Stanley (\cite{Stanley-hv}); 
see \cite{Bayer-Ehr} for formulas in terms of the flag vector.
Since it depends only on the flag vector, the toric $h$-vector of the
multiplex is the same as the toric $h$-vector of the 
pyramid over the appropriate polygon.
The $h$-vector of $M^{d,n}$ is thus
$$h(M^{d,n})=(1, 1, 1, \ldots, 1, 1) + (n-d)(0, 1, 1, \ldots, 1, 0).$$

A $d$-dimensional polytope $P$ is {\em elementary} if and only if the 
second and third entries in the toric $h$-vector are equal.  In terms
of the flag vector, this says
$f_{02}(P) - 3f_2(P) + f_1(P) -df_0(P)+{d+1 \choose 2} = 0$.
(See \cite{Kalai-NATO}.)

\begin{corollary}
For every $n\ge d\ge 2$, the multiplex $M^{d,n}$ is an elementary polytope.
\end{corollary}

Every face of a multiplex is a multiplex, but not all multiplexes are 
faces of higher dimensional multiplexes.

\begin{proposition}\label{multiplex-face}
The multiplex $M^{d,n}$ is a proper face of some multiplex if and only if
$d\le n\le 2d-1$.
\end{proposition}
\begin{proof}
Assume $M$ is a $d$-dimensional multiplex with $n+1$ vertices, and suppose
$M$ is a proper face of some higher dimensional multiplex $\widehat{M}$.
Then $M$ is a facet of a multiplex of dimension $d+1$, namely, any 
$(d+1)$-dimensional face of $\widehat{M}$ containing $M$.
According to the description of the facets in the definition of multiplex, 
$M$ has at most $2d$ vertices.
So $d+1\le n+1\le 2d$.

To prove the converse, let $Q$ be the $(d+1)$-dimensional multiplex with
$2d+1$ vertices.
The facets of $Q$ are 
$$F_i = \conv\{x_{i-d}, x_{i-d+1}, \ldots, x_{i-1}, x_{i+1}, x_{i+2}, \ldots,
x_{i+d}\},$$
for $0\le i\le 2d$.
For $1\le i\le d$, $F_i$ is a facet of $Q$ with $i+d$ vertices.
Thus among the proper faces of $Q$ are the $d$-multiplexes with $n+1$ vertices
for every $n$, $d\le n\le 2d-1$.
\end{proof}

Note that by Theorem~\ref{multi-backgr}, every polygon is a multiplex.
However, Proposition~\ref{multiplex-face} says that the two-dimensional faces of
higher dimensional multiplexes can only be triangles and quadrilaterals.

We turn now to the graphs of multiplexes.

\begin{theorem}[Bisztriczky \cite{Bisz-mult}]
Let $M^{d,n}$ be the multiplex with vertex set $\{x_0,x_1,\ldots, x_n\}$.
The edges of $M^{d,n}$ are 
\begin{itemize}
\item $\conv\{x_i,x_j\}$, where $0\le i < j \le n$ and $j-i\le d-2$
\item $\conv\{x_i,x_{i+d}\}$, where $0\le i  \le n-d$
\item $\conv\{x_0,x_{d-1}\}$
\item $\conv\{x_{n-d+1},x_n\}$.
\end{itemize}
\end{theorem}

There is a significant literature on the subject of reconstructing a polytope
from its graph.
For simplicial polytopes, the reconstruction problem can be phrased as the
problem of deciding when the vertices of a complete subgraph span a face of
the polytope.
The latter problem also arises when studying triangulations of polytopes.

\begin{theorem} \label{complete}
Every complete subgraph of the graph of a multiplex is the graph of a face
of the multiplex.
\end{theorem}
\begin{proof}
Let $G$ be the graph of the multiplex $M^{d,n}$.
Let $X=\{x_{\ell_0},x_{\ell_1},\ldots, x_{\ell_t}\}$ ${ } \subseteq\vertx(P)$ with
$\ell_0<\ell_1<\cdots<\ell_t$ such that for all $i,j\in\{0,1,\ldots,t\}$,
$\{x_{\ell_i},x_{\ell_j}\}$ is an edge of $G$.

{\em Case 1.}  If $\ell_0=0$, then $\ell_t\le d$. 
Note that $\{x_1,x_d\}$ is not an edge of $G$, so it is not a subset of $X$.
Thus $X$ is contained in either the facet $F_0$ or the facet $F_1$ of the
multiplex $M^{d,n}$.
Since both those facets are simplices, $\conv(X)$ is a simplex that is a face of
$M^{d,n}$.

If $\ell_t = n$, the argument is similar.

{\em Case 2.} Assume $0< \ell_0 < \ell_t < n$.
Note that each vertex $x_i$ of $M^{d,n}$ is contained in exactly the facets
$F_{i-d+1}$, $F_{i-d+2}$, \ldots, $F_{i-1}$, $F_{i+1}$, $F_{i+2}$, \ldots,
$F_{i+d-1}$.  
Here we use the conventions that $F_i=F_0$ if $i<0$, and $F_i=F_n$ if $i>n$.
(The self-duality of $M^{d,n}$ is expressed in the combinatorial description.)
Thus the set $X$ is contained in the facet $F_j$ if and only if 
$\ell_t-d+1\le j\le \ell_0+d-1$ and $j\not\in\{\ell_0,\ell_1,\ldots,\ell_t\}$.
(This statement is not true if $\ell_0=0$ or $\ell_t=n$.)
Let $a=\min\{0, \ell_t-d+1\}$ and $b=\max\{n,\ell_0+d-1\}$.
Then we can write the set $J=\{j: X\subset F_j\}$ as 
$J=[a,b]\setminus \{\ell_0,\ell_1,\ldots,\ell_t\}$.
The set of vertices in $\cap_{j\in J}F_j$ is exactly $X$, i.e., 
$\cap_{j\in J}F_j=\conv(X)$.
To check this, we need to know that $\{a,b\}\subseteq J$.
If $a=0$, then $a\in J$.
Otherwise note that since $x_{\ell_0}$ and $x_{\ell_t}$ are assumed to be
adjacent, either $\ell_t-d+1 < \ell_0$ (so $a\in J$), or $\ell_t = \ell_0+d$.
In the latter case, $a=\ell_t-d+1=\ell_0+1$. 
Since $x_{\ell_0+1}$ is not adjacent to 
$x_{\ell_t}=x_{\ell_0+d}$, $a=\ell_0+1\not\in \{\ell_0,\ell_1,\ldots,\ell_t\}$,
so $a\in J$.
Similarly, $b\in J$.
Now suppose $i\not\in \{\ell_0,\ell_1,\ldots,\ell_t\}$; we wish to show that
$x_i\not\in \cap_{j\in J}F_j$.
If $i<a$, then $x_i\not\in F_b$.
If $i>b$, then $x_i\not\in F_a$.
If $a\le i\le b$, then $x_i\in \cap_{j\in J}F_j$ if and only if
$i\in \{\ell_0,\ell_1,\ldots,\ell_t\}$.
Thus if $i\not\in \{\ell_0,\ell_1,\ldots,\ell_t\}$, then
$x_i\not\in \cap_{j\in J}F_j$, i.e., 
$\cap_{j\in J}F_j=\conv(X)$.

Now observe that this argument is also valid for every subsequence of 
$\ell_0<\ell_1<\cdots<\ell_t$.
Thus, for every $Y\subseteq X$, $\conv(Y)$ is a face of $M^{d,n}$, and
so also a face of $\conv(X)$.
Therefore $\conv(X)$ is a simplex.
\end{proof}

This last theorem enables us to count the number of faces that are
simplices.
\begin{proposition}
For $n>d$, the number of $i$-dimensional simplex faces of the multiplex 
$M^{d,n}$ is
$${d+1\choose i+1}-{d-3\choose i-3}+(n-d)\left[{d-1\choose i}-{d-3\choose i-2}
\right] .$$
\end{proposition}
\begin{proof}
Count the number of complete subgraphs with vertex set 
$X=\{x_{\ell_0},x_{\ell_1},\ldots, x_{\ell_i}\}$, with
$\ell_0<\ell_1<\cdots<\ell_i$.

$\bullet$ For $\ell_0=0$, the complete subgraphs are obtained by choosing
$i$-element sets from
$\{1,2,\ldots,d\}$ not containing the pair $\{1,d\}$, since $\{x_1,x_d\}$ is
the only nonedge using these indices.
There are ${d\choose i}-{d-2\choose i-2}$ such complete subgraphs.

$\bullet$ For $1\le \ell_0\le n-d-1$, the complete subgraphs are obtained 
by choosing $i$-elements sets from $\{\ell_0+1,\ell_0+2,\ldots, \ell_0+d-2,
\ell_0+d\}$ not containing the pair $\{\ell_0+1,\ell_0+d\}$.
There are ${d-1\choose i}-{d-3\choose i-2}$ such complete subgraphs for each
$\ell_0$, or 
$(n-d-1)[{d-1\choose i}-{d-3\choose i-2}]$ altogether.

$\bullet$ For $\ell_0=n-d$, the complete subgraphs are obtained by choosing
$i$-element sets from
$\{n-d+1,n-d+2,\ldots,n-2,n\}$.
There are ${d-1\choose i}$ such complete subgraphs.

$\bullet$ For $n-d+1\le \ell_0\le n-i$, the complete subgraphs are obtained by 
choosing $i$-element sets from
$\{\ell_0+1,\ell_0+2,\ldots,n-1,n\}$.
There are ${n-\ell_0\choose i}$ such complete subgraphs.

Altogether the number of $i$-faces that are simplices is thus
\begin{eqnarray*}
\lefteqn{{d\choose i}-{d-2\choose i-2}+(n-d-1)\left[{d-1\choose i}-{d-3\choose i-2}\right] + {d-1\choose i}+\sum_{j=i}^{d-1}{j\choose i}}\\
&=& 
{d\choose i}-{d-2\choose i-2}+(n-d)\left[{d-1\choose i}-{d-3\choose i-2}\right] + {d-3\choose i-2}+{d\choose i+1}\\
&=& 
{d+1\choose i+1}-{d-3\choose i-3}+(n-d)\left[{d-1\choose i}-{d-3\choose i-2}\right].
\end{eqnarray*}
\end{proof}

Thus, the number of nonsimplex $i$-faces of $M^{d,n}$ is 
${d-3\choose i-3}+(n-d){d-3\choose i-2}$.

Bisztriczky describes certain quadrilateral two-faces of the multiplex.

\begin{proposition}[\cite{Bisz-mult}] 
\label{quadril}
For $n>d\ge 3$, let $M^{d,n}$ be the multiplex with ordered list of vertices
$x_0$, $x_1$, \ldots, $x_n$.
Then for each $i$, \mbox{$0\le i\le n-d-1$}, $\conv\{x_i,x_{i+1},x_{i+d},x_{i+d+1}\}$
is a quadrilateral two face with diagonals \linebreak $\conv\{x_i,x_{i+d+1}\}$ and
$\conv\{x_{i+1},x_{i+d}\}$.
\end{proposition}

  From the comment above, we know that $M^{d,n}$ has exactly $n-d$ nonsimplex
two-faces.  Bisztriczky's proposition accounts for all of them.

The following propositions have analogues for ordinary polytopes, and
those are proved in section~\ref{sec-ord}.
The odd-dimensional multiplexes are also ordinary polytopes.
The proofs of Propositions~\ref{chromatic}--\ref{diameter} carry through
for the even-dimensional multiplexes as well.
The first is of special interest because Kalai (\cite{Kalai-NATO}) conjectures
that all elementary $d$-polytopes are $(d+1)$-colorable.

\begin{proposition}
The chromatic number of the graph of the multiplex $M^{d,n}$ is $d$, if $n>d$,
and $d+1$, if $n=d$ (in which case $M^{d,n}$ is the $d$-simplex).
\end{proposition}

\begin{proposition}
For every $n\ge d \ge 2$, the multiplex $M^{d,n}$ has a Hamiltonian cycle.
\end{proposition}

\begin{proposition}
For every $n\ge d \ge 2$, the multiplex $M^{d,n}$ has diameter 
$\lceil n/d\rceil$.
\end{proposition}

\section{Ordinary polytopes}\label{sec-ord}
Given an ordered set $V=\{x_0,x_1,\ldots,x_n\}$, a subset $Y\subseteq V$ is
called a {\em Gale subset} if between any two elements of $V\setminus Y$
there is an even number of elements of $Y$.
A polytope $P$ with ordered vertex set $V$ as above is a {\em Gale polytope}
if the set of vertices of each facet is a Gale subset.
\begin{definition}[\cite{Bisz}] {\em
An {\em ordinary polytope} is a Gale polytope such that each facet is a 
multiplex with the induced order on the vertices.
}\end{definition}
The definition of ordinary polytope is due to Bisztriczky.
His choice of the term ``ordinary'' stemmed from a belief that they arise
by choosing vertices on a convex ordinary space curve, but this is not
understood in dimensions higher than three.
See \cite{Bisz3,Bisz}.

The combinatorics of three-dimensional ordinary polytopes differs considerably
from that of higher dimensional ordinary polytopes. 
In this paper we consider only ordinary polytopes of dimension at least four.
Bisztriczky (\cite{Bisz}) defines these as above, and proves a number of
results on their combinatorics.
However, it is Dinh (\cite{Dinh}) who proves their existence in Euclidean
space.
We use the following theorems of Bisztriczky and Dinh.

\begin{theorem}[\cite{Bisz}]
Let $P$ be an ordinary $d$-polytope with ordered vertices 
$x_0$, $x_1$, \ldots, $x_n$.
Assume $n\ge d\ge 4$.  Then
\begin{enumerate}
\item If $d$ is even, then $P$ is cyclic.
\item If $d$ is odd, then there exists an integer $k$ ($d\le k\le n$) such that 
      the vertices sharing an edge with $x_0$ are exactly $x_1$, $x_2$, \ldots,
      $x_d$, and the vertices sharing an edge with $x_n$ are exactly
      $x_{n-1}$, $x_{n-2}$, \ldots, $x_{n-k}$.
\end{enumerate}
\end{theorem}
The integer $k$ guaranteed by this theorem is called the {\em characteristic}
of the ordinary polytope.

\begin{theorem}[\cite{Dinh}]
For every $n\ge k\ge d=2m+1\ge 5$, there exists an ordinary $d$-polytope with
$n+1$ vertices and characteristic $k$.
\end{theorem}

\begin{theorem}[\cite{Bisz}]\label{cyclic}
Given a triple of integers $(n,k,d)$ with $d$ odd and $n\ge k\ge d\ge 5$,
all ordinary $d$-polytopes with $n+1$ vertices and characteristic $k$ 
are combinatorially isomorphic.

If $P$ is an ordinary $d$-polytope with $n+1$ vertices and characteristic $k=n$,
then $P$ is a cyclic polytope.

If $P$ is an ordinary $d$-polytope with characteristic $k=d$, then $P$ is a
multiplex.
\end{theorem}

Write $P^{d,k,n}$ for the $d$-dimensional ordinary polytope with $n+1$ vertices
and characteristic $k$.
The appendix gives the face lattice of $P^{5,7,9}$; it was computed using
the software package {\em Polymake} (\cite{polymake}).

We do not have a formula for the complete flag vector of $P^{d,k,n}$.
However, Dinh (\cite{Dinh}) computes the $f$-vector of $P^{d,k,n}$, and
here we  compute $f_{02}$, and show that ordinary polytopes, other than
multiplexes, are not elementary.

We use the following description of the facets of $P^{d,k,n}$, due to Dinh.
\begin{theorem}[\cite{Dinh}]\label{ordfacet}
Let $n$, $k$, $d$ and $m$ be integers such that $n\ge k\ge d=2m+1\ge 5$.
Let $P^{d,k,n}$ be an ordinary $d$-polytope with characteristic $k$ and
ordered vertices $x_0$, $x_1$, \ldots, $x_n$.
Then the facets of $P^{d,k,n}$ are $\conv(X)$, where 
$$X=\{x_i,x_{i+1},\ldots,x_{i+2r-1}\}\cup Y\cup \{x_{i+k},x_{i+k+1},\ldots,
x_{i+k+2r-1}\},$$
where $i\in{\bf Z}$, $1\le r\le m$, $Y$ is a paired $(d-2r-1)$-element
subset of $\{x_{i+2r+1},x_{i+2r+2},\ldots, x_{i+k-2}\}$, and $|X|\ge d$.
Here a paired subset is one whose index set can be written as a disjoint union
of pairs of consecutive integers.
\end{theorem}

By definition each of these facets is a multiplex with the induced vertex
ordering.
Note that each facet (a $(d-1)$-dimensional multiplex) has at most
$d+2m-1=2d-2$ vertices.
Thus, Proposition~\ref{multiplex-face} generalizes to
\begin{quote}
{\em The multiplex $M^{d,n}$ is a proper face of some ordinary polytope if 
and only if $d\le n\le 2d-1$.}
\end{quote}
The facets of ordinary polytopes are small, simplex-like polytopes.
The duals of ordinary polytopes have facets with many vertices, however.
It can be shown that each facet of the dual of an ordinary $(2m+1)$-polytope
of characteristic $k$ has at least $3{k-m-3\choose m-1}$ vertices.
This is the number of those facets of $P^{d,k,n}$ that contain the vertex $x_1$
and fit the description of Theorem~\ref{ordfacet} with $r=1$.

The two-dimensional faces of $P^{d,k,n}$ are exactly the two-dimensional
faces of its facets.
 From Proposition~\ref{quadril} and the comment following it, we know all
the nontriangular two-faces.

\begin{proposition}\label{two-faces}
Let $n$, $k$, $d$ and $m$ be integers such that $n\ge k\ge d=2m+1\ge 5$.
Let $P^{d,k,n}$ be an ordinary $d$-polytope with characteristic $k$ and
ordered vertices $x_0$, $x_1$, \ldots, $x_n$.
The two-dimensional faces of $P^{d,k,n}$ that are not triangles are 
exactly the sets $\conv\{x_i,x_{i+1},x_{i+k},x_{i+k+1}\}$, for 
$0\le i\le n-k-1$.
\end{proposition}
\begin{proof}
Dinh (\cite{Dinh}) proves that these are indeed two-faces of $P^{d,k,n}$,
but does not show that they are the only nontriangular two-faces.
Consider a facet $F$ as given in Theorem~\ref{ordfacet}, and apply 
Proposition~\ref{quadril} (and the subsequent comment) to give all its
nontriangular two-faces.
Renumber the vertices of the facet $F$ as $z_0$, $z_1$, \ldots.
Thus the $d-2r-1$ elements of $Y$ are numbered $z_{2r}$ through $z_{d-2}$,
and $\{x_{i+k},x_{i+k+1},\ldots, x_{i+k+2r-1}\}=
\{z_{d-1},z_d,\ldots, z_{d+2r-2}\}$.
The quadrilateral faces in $F$ are of the form
$\conv\{z_j,z_{j+1},z_{j+d-1},z_{j+d}\}$.
Such a quadrilateral contains no element of $Y$, and is of the form \linebreak
$\conv\{x_\ell,x_{\ell+1},x_{\ell+k},x_{\ell+k+1}\}$, with $i\le\ell\le i+2r-2$.
Considering all the facets of $P^{d,k,n}$, $\ell$ can range from 0 to $n-k-1$.
\end{proof}

Dinh also computes the $f$-vectors of ordinary polytopes.
At the moment we need only $f_1$.
\begin{proposition}[\cite{Dinh}]\label{f_1}
Let $n$, $k$, $d$ and $m$ be integers such that $n\ge k\ge d=2m+1\ge 5$.
Then $f_1(P^{d,k,n})={k+1\choose 2}+(n-k)(k-1)$.
\end{proposition}

\begin{proposition}
For every $n\ge k\ge  d =2m+1\ge 5$,
the ordinary polytope $P^{d,k,n}$ is elementary if and only if $k=d$, in 
which case the ordinary polytope is itself a multiplex.
\end{proposition}
\begin{proof}
Proposition~\ref{two-faces} implies that for the ordinary polytope $P^{d,k,n}$,
\linebreak
$f_{02}(P^{d,k,n})=3f_2(P^{d,k,n})+(n-k)$.
Write $$\beta=f_{02}(P^{d,k,n})-3f_2(P^{d,k,n})+f_1(P^{d,k,n})-df_0(P^{d,k,n})+
{d+1\choose 2}.$$
Then 
$\beta= (n-k) + {k+1\choose 2} + (n-k)(k-1) - d(n+1) +{d+1\choose 2}$.
The ordinary polytope $P^{d,k,n}$ is elementary if and only if $\beta = 0$.
If $k\ne d$, then solving $\beta = 0$ for $n$ gives $n=(k+d-1)/2$.
But $(k+d-1)/2 \le (2k-1)/2 < k \le n$, so this is not possible.
Thus $\beta=0$ implies $k=d$, i.e., the ordinary polytope is a multiplex.
We have already seen that every multiplex is elementary.
\end{proof}

By Dinh's construction (\cite{Dinh}) of ordinary polytopes, they can be
realized as rational polytopes.
It would be interesting to compute the toric $h$-vector of ordinary polytopes.
The quantities above are enough for the case of dimension five.
\begin{theorem}
Let $n\ge k\ge 5$.  The toric $h$-vector of $P^{5,k,n}$ is
$$(1, n-4, {n-3\choose 2}-{n-k+1\choose 2}, {n-3\choose 2}-{n-k+1\choose 2}, 
n-4,1).$$
\end{theorem}
Among all 5-polytopes with $h_1=n-4$ (that is, $f_0=n+1$), the smallest 
possible $h_2$ is $n-4$, and this is achieved by the multiplex 
$M^{5,n}=P^{5,5,n}$.
Among all 5-polytopes with $h_1=n-4$, the largest possible $h_2$ is
${n-3\choose 2}$, and this is achieved by the cyclic polytope $P^{5,n,n}$.
Thus, the toric $h$-vectors of $P^{5,k,n}$ are nicely distributed between
the extreme toric $h$-vectors having $h_1=n-4$.
The formula for $h_2$ generalizes for ordinary polytopes of higher (odd)
dimension: $h_2(P^{d,k,n})={n-d+2\choose 2}-{n-k+1\choose 2}$; note that
${n-d+2\choose 2}$ is $h_2$ for the cyclic $d$-polytope with $n+1$ vertices.

The description of the graph of a multiplex extends naturally to ordinary 
polytopes.

\begin{theorem}
For $n\ge k\ge d=2m+1\ge 5$,
let $P^{d,k,n}$ be the ordinary polytope with vertex set 
$\{x_0,x_1,\ldots, x_n\}$.
The edges of $P^{d,k,n}$ are 
\begin{itemize}
\item $\conv\{x_i,x_j\}$, where $0\le i < j \le n$ and $j-i\le k-2$
\item $\conv\{x_i,x_{i+k}\}$, where $0\le i  \le n-k$
\item $\conv\{x_0,x_{k-1}\}$
\item $\conv\{x_{n-k+1},x_n\}$.
\end{itemize}
\end{theorem}
\begin{proof}
We consider the pairs not listed in the statement of the theorem.
These fall into two categories.  First are those pairs $\{x_i,x_j\}$ with
$j-i\ge k+1$.
 From Theorem~\ref{ordfacet} every facet containing $x_i$ and $x_j$, with
$j-i\ge k+1$, also contains the nonempty set of vertices
$\{x_{i+k},x_{i+k+1},\ldots, x_{j-1}\}$.
Thus, $\conv\{x_i,x_j\}$ is not a face (edge) of $P^{d,k,n}$.
The number of these pairs is $\sum_{i=0}^{n-k-1} (n-k-i) = 
\sum_{\ell=1}^{n-k} \ell = {n-k+1\choose 2}$.
The other pairs not listed are $\{x_i,x_j\}$ with $j-i=k-1$, $i\ne 0$,
and $j\ne n$.
For such pairs $\conv\{x_i,x_j\}$ is not an edge of $P^{d,k,n}$, because
it is a diagonal of a two-dimensional face as described in 
Proposition~\ref{two-faces}.
The number of these pairs is $n-k$.
The number of pairs listed in the statement of the theorem is thus
${n+1\choose 2}-{n-k+1\choose 2}-(n-k)$ ${}={k+1\choose 2}+(n-k)(k-1)$,
which is the number of edges of $P^{d,k,n}$, as computed by Dinh
(Proposition~\ref{f_1}).
So all the listed pairs are edges of $P^{d,k,n}$.
\end{proof}

Note that the graphs of ordinary polytopes are dimensionally 
ambiguous.
For every odd $d$ and $d'$ between 5 and $k$, the graphs of $P^{d,k,n}$ and
$P^{d',k,n}$ are isomorphic.

\begin{proposition}\label{chromatic}
For every $n\ge k\ge  d =2m+1\ge 5$,
the chromatic number of the graph of the ordinary polytope $P^{d,k,n}$ is $k$,
if $n>k$, and $k+1$, if $n=k$ (in which case $P^{d,k,n}$ is a cyclic polytope).
\end{proposition}
\begin{proof}
If $n=k$, then $P=P^{d,k,n}$ is a cyclic polytope, and the graph of $P$
is the complete graph on $k+1$ vertices, so its chromatic number is $k+1$.
So assume $n>k$.
Let $G$ be the graph of $P$, with vertex set $\{x_0,x_1,\ldots, x_n\}$.
Assign colors from the set $\{0,1,\ldots,k-1\}$ to the vertices of $G$ as
follows:
$$\lambda(x_i)=\left\{ \begin{array}{cl}
                       k-1 & \mbox{if $i=0$ or $i=n$} \\
		       i \bmod k-1 & \mbox{if $1\le i\le n-1$} \end{array}
		       \right. $$
Since $x_0$ and $x_n$ are not adjacent in $G$, every edge containing $x_0$
or $x_n$ is assigned two different colors.
If $i<j$ and $j-i\le k-2$, or $j-i=k$, then
$j \not\equiv i \pmod{k-1}$, so $\lambda(x_j)\ne \lambda(x_i)$.
Thus adjacent vertices have different colors, so $\lambda$ gives a
proper $k$-coloring of $G$.
Now $G$ contains a complete subgraph on the vertex set 
$\{x_0, x_1, \ldots, x_{k-1}\}$, so the chromatic number of $G$ is $k$.
\end{proof}

\begin{proposition}
For every $n\ge k\ge  d =2m+1\ge 5$, the ordinary polytope $P^{d,k,n}$ has a 
Hamiltonian cycle.
\end{proposition}
\begin{proof}
If $n$ is odd, the vertex sequence, $x_0$, $x_2$, $x_4$, \ldots, $x_{n-1}$,
$x_n$, $x_{n-2}$, $x_{n-4}$, \ldots, $x_3$, $x_1$, $x_0$, gives a Hamiltonian
cycle.
If $n$ is even, the vertex sequence, $x_0$, $x_2$, $x_4$, \ldots, $x_{n-2}$,
$x_n$, $x_{n-1}$, $x_{n-3}$, \ldots, $x_3$, $x_1$, $x_0$, gives a Hamiltonian
cycle.
\end{proof}

\begin{proposition}\label{diameter}
For every $n\ge k\ge d=2m+1 \ge 5$, the ordinary polytope $P^{d,k,n}$ has 
diameter $\lceil n/k\rceil$.
\end{proposition}
\begin{proof}
For $i<j$, the vertices $x_i$ and $x_j$ are adjacent if $j-i\le k-2$
or if $j-i=k$.
In addition $x_0$ and $x_{k-1}$ are adjacent, and $x_n$ and $x_{n-k+1}$ are
adjacent.
So usually $x_i$, $x_{i+k}$, $x_{i+2k}$, \ldots, $x_{i+tk}$, $x_j$, with
$t=\lfloor j-i-1/k\rfloor$, gives an $x_i$--$x_j$ path of length
$t+1\le \lceil n/k\rceil$.
This is valid as long as $j-i\not\equiv -1 \pmod{k}$.
If $j-i\equiv -1\pmod{k}$, then
$x_i$, $x_{i+k}$, $x_{i+2k}$, \ldots, $x_{i+tk}$, $x_{i+tk+1}$, $x_j$, with
$t=\lfloor(j-i)/k\rfloor=(j-i-k+1)/k$, gives a path of length $t+2$.
If $j-i<n-1$, then $t+2=(j-i+k+1)/k\le \lceil n/k\rceil$.
The remaining cases are the paths
$x_0$, $x_{k-1}$, $x_{2k-1}$, \ldots, $x_{n-1}$ and
$x_1$, $x_{k+1}$, $x_{2k+1}$, \ldots, $x_{n+1-k}$, $x_n$, if 
$n\equiv 0 \pmod{k}$, and 
$x_0$, $x_k$, $x_{2k}$, \ldots, $x_{tk}$, $x_n$, with $t=\lfloor n/k\rfloor$, if
$n\equiv -1 \pmod{k}$.
These are all of length $\lceil n/k\rceil$.
Clearly, the $x_0$--$x_n$ path given is the shortest $x_0$--$x_n$ path.
So the diameter is exactly $\lceil n/k\rceil$.
\end{proof}

By Theorem~\ref{cyclic} the class of ordinary polytopes includes the 
cyclic polytopes of odd dimension.  
Cyclic polytopes have played an important role in the combinatorial
study of simplicial polytopes (e.g., in \cite{McM}), and more recently, 
in the study of triangulations of polytopes (e.g., in \cite{rambau-bauesI}).
Cyclic polytopes are {\em neighborly}, that is, every 
$\lfloor d/2\rfloor$-element set of vertices is the vertex set of a face of
the polytope.
In particular, the graph of a neighborly $d$-polytope for $d\ge 4$ is the 
complete graph.
Thus no ordinary polytopes other than the cyclic polytopes are neighborly.
A generalization of neighborliness is studied in \cite{Bayer-weakly}.
A polytope is {\em weakly neighborly} if every set of $k+1$ vertices
is contained in a face of dimension at most $2k$, for all $k$.
It is natural to ask, then, if ordinary polytopes are weakly neighborly.
The answer is no, almost always.
If $n\ge k+2$, then $x_0$ and $x_n$ are not on a common two-face of the
ordinary polytope $P^{d,k,n}$.
For $n=k+1$, $k>d\ge 5$ ($d$ odd), $P^{d,k,n}$ is not weakly neighborly;
for example, $\{x_1,x_3,x_5\}$ is not contained in a facet of $P^{5,6,7}$.
If $n=d+2$, then $x_0$ and $x_n$ are not on a common two-face of the 
multiplex $M^{d,n}$.
The only weakly neighborly ordinary polytopes are the cyclic
polytopes and the multiplexes $M^{d,d+1}$, which are $(d-2)$-fold pyramids
over quadrilaterals.
Theorem~\ref{complete} says that every complete subgraph of the graph of a 
multiplex is the graph of a face of the multiplex.
This fails in general for ordinary polytopes.

We turn now to the $f$-vectors of ordinary polytopes.
These are computed by Dinh.
\begin{theorem}[\cite{Dinh}]
Let $n\ge k\ge d=2m+1\ge 5$.
The number of \linebreak $i$-dimensional faces of the ordinary polytope $P^{d,k,n}$ is
$$f_i(P^{d,k,n})=\phi_i(d,k)+(n-k)c_i(d,k),$$
where $\phi_i(d,k)$ is the number of $i$-faces of the cyclic $d$-polytope with
$k+1$ vertices,
$$\phi_i=\left\{ \begin{array}{cl}
	  {k+1\choose i+1} & \mbox{for $0\le i\le m-1$} \\
          \sum_{j=0}^m \left[ {j\choose d-1-i} + {d-j\choose d-1-i}\right]
	  {k-d+j\choose j} & \mbox{for $m\le i\le d-1$} \end{array}
	  \right. ,$$
and $c_i(d,k)=f_i(P^{d,k,n+1})-f_i(P^{d,k,n})$ is given by 
\begin{enumerate}
\item $\displaystyle c_i={k-1\choose i}$, for $1\le i<m$,
\item $\displaystyle c_m={k-1 \choose m}-{k-2-m\choose m}$,
\item $\displaystyle c_i=\sum_{j=i-m}^{\lfloor i/2\rfloor}
      (2N(k-1,j,i)-N(k-2,j,i))-
      \sum_{j=i-m}^{\lfloor (i-1)/2\rfloor} N(k-3,j,i-1)$ 
      $\displaystyle {} -
      \sum_{j=i-m-1}^{\lfloor (i-2)/2\rfloor} N(k-3,j,i-2)-
      \sum_{r=0}^{i-m} N(k-3-2r,i-m-r,i-2r)$, \\  \\
      for $m<i<2m=d-1$, and

\item $\displaystyle c_{d-1}=c_{2m}={k-2-m\choose m-1}$.
\end{enumerate}
Here 
$$N(s,t,u)={u-t\choose t}{s-u+t\choose u-t}+{u-1-t\choose t}{s-u+t\choose u-1-t}.$$
\end{theorem}

Thus, for fixed $d$ and $k$, the $f$-vectors of the ordinary polytopes 
$P^{d,k,n}$ lie on a line.

\begin{conjecture}
Let $d$ be an odd integer, $d\ge 5$.
The set of $f$-vectors of all ordinary $d$-polytopes spans the Euler 
hyperplane (given by $\sum_{i=0}^{d-1}f_i=2$).
A spanning set consists of the ordinary polytopes $P^{d,d+\lfloor i/2\rfloor,
d+i}$, for $1\le i\le d$.
\end{conjecture}
The conjecture has been verified on computer for odd $d$, $5\le d\le 37$.

The flag vectors of ordinary polytopes satisfy many linear relations that do 
not hold for all polytopes, however.
In particular, the self-duality of multiplexes gives the following equalities.
Let $S=\{s_1,s_2,\ldots, s_{r-1},s_r\}$, and 
$S'=\{s_r-1-s_1,s_r-1-s_2,\ldots, s_r-1-s_{r-1},s_r\}$.
Then for every ordinary polytope $P$ of odd dimension greater than $s_r$,
$f_S(P)=f_{S'}(P)$.
For example, for ordinary 5-dimensional polytopes, $f_{03}=f_{23}$, a relation
that fails for arbitrary 5-polytopes.
In fact, in dimension five, the flag vectors of ordinary polytopes depend 
linearly on the $f$-vectors.

\pagebreak

\renewcommand{\thesection} {}

\section{Appendix}

Here is the face lattice of the ordinary polytope $P^{5,7,9}$.
Faces are listed by their sets of vertices, from the vertex set
$\{0, 1, \ldots, 9\}$.
The $f$-vector of this polytope is $f(P^{5,7,9})=(10, 40, 76, 70, 26)$.

{\addtolength{\parskip}{3pt}

\vspace*{12pt}

\noindent Facets:

01234, 01245, 01256, 02345, 02356, 02367, 03456, 03467, 04567, 23459, 

23569, 23679, 34569, 34679, 34789, 45679, 45789, 56789, 013478, 

014578, 015678, 123489, 124589, 125689, 0123789, 0126789

\vspace*{6pt}

\noindent 3-dimensional faces:

0123, 0124, 0125, 0126, 0134, 0145, 0156, 0234, 0235, 0236, 0237, 0245,

0256, 0267, 0345, 0346, 0347, 0356, 0367, 0456, 0457, 0467, 0567, 1234,

1245, 1256, 1348, 1458, 1568, 2345, 2349, 2356, 2359, 2367, 2369, 2379,

2459, 2569, 2679, 3456, 3459, 3467, 3469, 3478, 3479, 3489, 3569, 3679,

3789, 4567, 4569, 4578, 4579, 4589, 4679, 4789, 5678, 5679, 5689, 5789,

6789, 01378, 01478, 01578, 01678, 12389, 12489, 12589, 12689, 012789

\vspace*{6pt}

\noindent 2-dimensional faces:

012, 013, 014, 015, 016, 023, 024, 025, 026, 027, 034, 035, 036, 037,

045, 046, 047, 056, 057, 067, 123, 124, 125, 126, 134, 138, 145, 148,

156, 158, 168, 234, 235, 236, 237, 239, 245, 249, 256, 259, 267, 269,

279, 345, 346, 347, 348, 349, 356, 359, 367, 369, 378, 379, 389, 456,

457, 458, 459, 467, 469, 478, 479, 489, 567, 568, 569, 578, 579, 589,

678, 679, 689, 789, 0178, 1289

\vspace*{6pt}

\noindent Edges:

01, 02, 03, 04, 05, 06, 07, 12, 13, 14, 15, 16, 18, 23, 24, 25, 26, 27,

29, 34, 35, 36, 37, 38, 39, 45, 46, 47, 48, 49, 56, 57, 58, 59, 67, 68,

69, 78, 79, 89}

\end{document}